# On the structure of a family of probability generating functions induced by shock models

## Satrajit Roychoudhury[*1] and Manish C. Bhattacharjee[2]

*New Jersey Institute of Technology*

**Abstract:** We explore conditions for a class of functions defined via an integral representation to be a probability generating function of some positive integer valued random variable. Interest in and research on this question is motivated by an apparently surprising connection between a family of classic shock models due to Esary et. al. (1973) and the negatively aging nonparametric notion of "strongly decreasing failure rate" (SDFR) introduced by Bhattacharjee (2005). A counterexample shows that there exist probability generating functions with our integral representation which are not discrete SDFR, but when used as shock resistance probabilities can give rise to a SDFR survival distribution in continuous time.

## 1. The Problem

A surprising connection between a family of classic shock models (Esary et al. [2]) and the negatively aging nonparametric notion of *strong decreasing failure rate* (SDFR) (Bhattacharjee [1]) is exploited in this paper to investigate necessary and sufficient conditions for a class of functions defined via an integral representation to be a probability generating function (p.g.f) of some positive integer valued random variable (r.v.). The question we investigate can be simply posed as follows. If $Q$ is a probability measure on the half line, under what conditions is the function defined by

$$(1.1) \qquad \int_0^\infty \frac{yz}{1 - z + yz} \, Q(dy), \qquad 0 < z < 1,$$

a p.g.f. of some positive integer valued r.v. $N$ ? This problem arose in the context of investigating a strong negatively aging property, introduced and studied by one of us (Bhattacharjee [1]), and its preservation by a class of shock models. We show that a weak converse to the preservation result also holds and provide a counter example to illustrate the breakdown of the full converse. Our work can thus be viewed as a non-traditional application of reliability ideas in the sense that the question posed can be conceived independently of the reliability theoretic framework, while ideas from the latter proves fruitful for the solution. The only other work we are aware of that is analogous in spirit, though entirely different from our work presented here, in terms of the specific question posed and results obtained, is a paper by Pakes [6] on the structure of a class of p.g.f.s and its connections to Markov branching processes.

---

*Satrajit Roychoudhury is currently in Schering Plough Research Institute.

[1]Schering Plough Research, Kenilworth, NJ 07033, USA, e-mail: satrajit.roychoudhury@spcorp.com

[2]Department of Mathematical sciences, New Jersey Institute of Technology, Newark, NJ 07102, USA, e-mail: manish.bhattacharjee@njit.edu

For any $y$ in $(0, 1)$, recognizing the integrand to be the p.g.f. of a geometric distribution over the positive integers, the answer is clearly affirmative if the support of the mixing distribution $\boldsymbol{Q}$ is no larger than $(0, 1]$. The case $y = 1$ corresponds to a mixing distribution degenerate at 1. For $y \in (0, 1]$, we can think of $N$ as conditionally geometric given $y$, so that (1.1) is the unconditional p.g.f of $N$, when the parameter $y$ is randomized over $(0, 1]$. In other words, if $\boldsymbol{Q}(0, 1] = 1$, then the function defined by (1.1) is a Bayesian's view of the p.g.f of $N$ when $y$ has a prior $\boldsymbol{Q}$. In fact if $\{X_1, X_2, \ldots\}$ is a sequence of binary exchangeable random variables, then for any integers $n$, $k$ such that $1 \le k \le n$, any $\{i_1, i_2, \ldots, i_k\} \subset \{1, 2, \ldots, n\}$, by de Finetti's classic result,

$$P\{X_{i_1} = X_{i_2} = \cdots = X_{i_k} = 1, \ X_j = 0 \text{ for } j \in \{1, \ldots, n\} \setminus \{i_1, \ldots, i_k\}\}$$
$$= \int_0^1 y^k (1 - y)^{n-k} \, \boldsymbol{Q}(dy),$$

for some probability measure $\boldsymbol{Q}$ in $(0, 1]$; it easily follows that the random variable

$$N = \inf\{n \ge 1 \ : \ X_n = 1\}$$

indeed has the p.g.f.

$$\int_0^1 \frac{yz}{1 - z + yz} \, \boldsymbol{Q}(dy)$$

for some unique probability measure $\boldsymbol{Q}$ supported by the unit interval.

However, if the support of $\boldsymbol{Q}$ on $\Re_+ \equiv [0, \infty)$ extends beyond $[0, 1]$, then the answer to our question is not obvious. Hence, we ask if $\boldsymbol{Q}(1, \infty) > 0$, can (1.1) still be a p.g.f. of a random variable $N$? The next section explores and answers this question.

## 2. Motivation and main results

Our problem is motivated by the following observations and results. Let

$$\bar{S}(t) = \sum_{k=0}^{\infty} \bar{P}_k e^{-\lambda t} \frac{(\lambda t)^k}{k!} \tag{2.1}$$

be the survival probability of the standard Esary-Marshall-Proschan (EMP) shock model (Esary et al. [2]), where failure is caused by shocks arising over time according to a homogeneous Poisson process $\{N(t); t \ge 0\}$ with intensity $\lambda > 0$ and the number of shocks($J$) to failure has tail (*shock resistance*) probabilities $\bar{P}_k = P(J > k)$, $k = 0, 1, 2, \ldots$. Esary et al. [2] proved that all of the standard nonparametric positive and negative aging properties of $J$ in discrete time are preserved by the survival distribution $S \equiv 1 - \bar{S}$ in continuous time.

**Definition 2.1.** (i) A non-discrete lifetime $X$ with distribution function $F$ has the SDFR(Strongly Decreasing Failure Rate) property if the tail (reliability) function $\bar{F}(t) = P(X > t)$ is a *completely monotone function* (Feller [3, 4]) on $[0, \infty)$.

(ii) A discrete non-negative integer valued random variable $X$ is *discrete* SDFR if its tail probability $\{u_n = P(X > n): \ n = 0, 1, 2, \ldots\}$, with $u_0 = 1$ is a *completely monotone sequence* (Feller [3, 4]).

A positive integer valued random variable $X$ has the discrete SDFR property if the non-negative random variable $(X - 1)$ is *discrete* SDFR.

**Lemma 2.1.** *The Laplace-Stieltjes transform of $S$ in (2.1) is given by*

$$L(s) = \phi\left(\lambda/(\lambda + s)\right), \quad s > 0,$$

*where $\phi(\cdot)$ is the p.g.f. of the random number $J$ of shocks to failure.*

*Proof.* Follows by routine computations; *viz.*, for $s > 0$,

$$
\begin{aligned}
s^{-1}\{1 - L(s)\} = \int_0^\infty e^{-st}\bar{S}(t)\,dt &= \sum_{k=0}^\infty \bar{P}_k \frac{\lambda^k}{(\lambda + s)^{k+1}}, \quad \text{using (2.1)} \\
&= \frac{1}{\lambda + s}\sum_{k=0}^\infty \bar{P}_k\left(\frac{\lambda}{\lambda + s}\right)^k.
\end{aligned}
$$

Now use the standard identity,

$$\sum_{k=0}^\infty \bar{P}_k z^k = \frac{1 - \phi(z)}{1 - z} \ \text{ with } \ z = \frac{\lambda}{\lambda + s} \in (0, 1). \qquad \square$$

The following characterization of the *non-discrete* SDFR property is known from Bhattacharjee [1].

**Theorem 2.1.** *$X \sim F$ is SDFR iff it has a representation*

$$F(t) = \int_0^\infty (1 - e^{-\lambda t})G(d\lambda)$$

*for some unique mixing distribution $G$ continuous at zero or, equivalently*

$$X \text{ is SDFR} \Leftrightarrow X \overset{d}{=} \frac{Y}{Z},$$

*where $\overset{d}{=}$ denotes equality in distribution, $Y$ and $Z$ are independent, $Y$ is exponentially distributed with unit mean and $Z$ is a positive random variable.*

*Sketch of proof.* By Bernstein's representation theorem of completely monotone functions (Feller [4]), there exists a unique, possibly substochastic, measure $G$ on $[0, \infty)$ such that for $t > 0$,

$$\bar{F}(t) = \int_0^\infty e^{-\lambda t}G(d\lambda) = G\{0\} + \int_{(0,\infty)} e^{-\lambda t}G(d\lambda).$$

We must have $G\{0\} = 0$, since $\bar{F}(t) \to 0$ as $t \to 0$. $\qquad \square$

If the shock resistance probabilities $\bar{P}_k$ in (2.1) are discrete SDFR, does the shock model survival probability inherit the SDFR property in continuous time? The answer is affirmative.

**Theorem 2.2.** *If $\bar{P}_k$ is SDFR, then $S$ is SDFR.*

*Proof.* Using Hausdorff's result (Feller [4], see Theorem 2, pages 222–224) on the equivalence of moment sequences of distributions on the unit interval $[0, 1]$ and the complete monotonicity property of such sequences,

$$\bar{P}_k \text{ is SDFR with } \bar{P}_0 = 1 \iff \bar{P}_k = \int_0^1 p^k \, \boldsymbol{F}(dp),$$

for some probability measure $\boldsymbol{F}$ on the unit interval. Accordingly from (2.1), we have,

$$\begin{aligned}
\bar{S}(t) &= \sum_{k=0}^{\infty} (\int_0^1 p^k \, \boldsymbol{F}(dp)) \, e^{-\lambda t} \frac{(\lambda t)^k}{k!} \\
&= \int_0^1 e^{-\lambda t(1-p)} \, \boldsymbol{F}(dp) = \int_0^{\lambda} e^{-\theta t} \boldsymbol{H}(d\theta),
\end{aligned}$$

for some probability measure $\boldsymbol{H}$. Note, $\lambda > 0$ implies that the support of $\boldsymbol{H}$ is contained in $[0, \infty)$, even though $\boldsymbol{F}$ has no mass outside the unit interval. Theorem 2.1 now implies that $S$ is SDFR. □

A necessary and sufficient condition for the EMP shock model survival distribution $S$ to be SDFR is given by the next result, that sets in perspective our question posed in Section 1.

**Theorem 2.3.** *The EMP Shock Model distribution function $S$ is SDFR iff the number of shocks to failure has a probability generating function $\phi(z) = Ez^J$ with a unique representation*

$$(2.2) \qquad \phi(z) = \int_0^{\infty} \frac{zy}{(1-z) + zy} \, \boldsymbol{Q}(dy)$$

*for some mixing distribution $\boldsymbol{Q}$ with support in $(0, \infty)$.*

*Proof.* By Lemma 2.1, for $s > 0$,

$$\begin{aligned}
\phi\left(\frac{\lambda}{\lambda + s}\right) &= E(e^{-sT}), \quad \text{where } T \text{ has tail } \bar{S}(t) \text{ in (2.1)} \\
&\equiv \int_0^{\infty} e^{-st} S(dt) \\
&= \int_0^{\infty} \int_0^{\infty} \theta e^{-(\theta + s)t} \, G(d\theta) \, dt \quad \text{(using Theorem 2.1)} \\
&= \int_0^{\infty} \frac{\theta}{\theta + s} \, G(d\theta) \quad \text{(using Fubini's Theorem)}.
\end{aligned}$$

Setting $z = \lambda/(\lambda + s) \in (0, 1)$ as $s \in (0, \infty)$, $\lambda > 0$; this yields

$$\begin{aligned}
\phi(z) &= \int_0^{\infty} \frac{\theta}{\theta + \lambda(\frac{1}{z} - 1)} \, G(d\theta) \\
&= \int_0^{\infty} \frac{\theta \lambda^{-1}}{\theta \lambda^{-1} + (\frac{1}{z} - 1)} \, G(d\theta), \quad \text{since } G\{0\} = 0 \\
&= \int_0^{\infty} \frac{zy}{(1-z) + zy} \, \boldsymbol{Q}(dy) \\
&= E\left(\frac{zY}{1 - z + zY}\right),
\end{aligned}$$

where $\Theta \sim G$ and $\frac{\Theta}{\lambda} \equiv Y \sim \boldsymbol{Q}$.

Note that, since $G$ is continuous at zero by Theorem 2.1, so must be $\boldsymbol{Q}$. Hence $\boldsymbol{Q}(0) = 0$. Also, since $\lambda > 0$, clearly the mixing distribution $\boldsymbol{Q}$ of the mixing variable $Y$ must have support contained in $(0, \infty)$. $\qquad\square$

At this point an obvious question is: what more can we say about $\boldsymbol{Q}$? In particular, what should be the support of $\boldsymbol{Q}$? Are there necessary and sufficient conditions on $\boldsymbol{Q}$ such that the right hand side of (1.1) is always a *probability generating function* (p.g.f.)? To explore such conditions, a nontrivial observation is that $\boldsymbol{Q}$ cannot have zero mass in $(0, 1]$. To check this claim, rewrite (2.2) as

$$\phi(z) = \int_0^\infty H(y, z) \boldsymbol{Q}(dy),$$

where the integrand

$$H(y, z) = \frac{zy}{1 - z + zy}, \quad y \geq 0, \ 0 < z \leq 1,$$

is concave $\uparrow$ in $y \geq 0$, for each $z \in (0, 1]$. Thus, since

$$\inf_{0 < y \leq 1} H(y, z) = H(0+, z) = 0, \text{ and } \inf_{y > 1} H(y, z) = H(1, z) = z,$$

we have,

$$\phi(z) \geq \int_1^\infty \inf_{y > 1} H(y, z) \boldsymbol{Q}(dy) = z \boldsymbol{Q}(1, \infty),$$

so that, $\boldsymbol{Q}(1, \infty) \leq z^{-1} \phi(z)$, for $0 < z \leq 1$; whence, it follows that

$$\boldsymbol{Q}(1, \infty) \leq \lim_{z \to 0+} \frac{\phi(z)}{z} = \phi'(0) = p_1 \equiv P(J = 1),$$

where $J$ is the number of shocks to failure in the EMP shock model (2.1), with tail probabilities $\bar{P}_k$. This must imply $\boldsymbol{Q}(0, 1] > 0$; for *if not*, then $1 \geq p_1 \geq \boldsymbol{Q}(1, \infty) = 1$. Thus $J = 1$ w.p. 1 and hence $\phi(z) = E(z^J) = z$, which in turn, in virtue of (2.2) implies $\boldsymbol{Q}\{1\} = 1$ and contradicts the hypothesis $\boldsymbol{Q}(0, 1] = 0$.

Similarly, using the monotonicity in $y$ of the integrand in (2.2), and by noting $\sup_{0 < y \leq 1} H(y, z) = H(1, z) = z$, and $\sup_{y > 1} H(y, z) = \lim_{y \to \infty} H(y, z) = 1$, we see that $\phi(z)$ defined in (2.2) satisfies

$$\phi(z) \leq z \boldsymbol{Q}(0, 1] + \boldsymbol{Q}(1, \infty).$$

Hence if the entire mass of $\boldsymbol{Q}$ lies in the unit interval $(0, 1]$, then $\phi(z) \leq z$, so that the graph of the p.g.f. $\phi(z)$ lies entirely below the diagonal in the unit square. In other words, the number $(J)$ of shocks to failure dominates the degenerate random variable constant at 1 in the *generating function order* (Müller and Stoyan [5]).

These observations still beg the question: can $\boldsymbol{Q}$ allocate positive mass to $(1, \infty)$, while $\phi(z)$ in (2.2) still remains a p.g.f.? In attempting to answer this question, the following lemma will prove useful.

**Lemma 2.2.** *For any non-negative integer valued random variable $N$ with distribution, $q_n \equiv P(N = n)$; $n = 0, 1, 2, \cdots$ we have*

$$E(1 - z)^N = \sum_{k=0}^\infty c_k z^k, \text{ where } c_k = (-1)^k \sum_{n=k}^\infty \binom{n}{k} q_n, \ 0 < z < 1.$$

*Proof.* Simply note,

$$E(1-z)^N = \sum_{n=0}^{\infty}(1-z)^n q_n \;=\; \sum_{n=0}^{\infty}\{\sum_{k=0}^{n}\binom{n}{k}(-z)^k\}q_n$$

$$=\; \sum_{k=0}^{\infty}\{\sum_{n=k}^{\infty}\binom{n}{k}q_n\}(-z)^k. \;\square$$

This leads to the following necessary condition.

**Theorem 2.4.** *If the integral in (1.1) is a probability generating function of a positive integer valued random variable $N$, then $\boldsymbol{Q}$ cannot have support beyond $(0,2)$.*

*Proof.* By Theorem 2.3, if the integral in (1.1) is to be a p.g.f., then it must correspond to the p.g.f. $\phi(z)$ of the number of shocks to failure in the EMP shock model (2.1). Hence,

$$\phi(z) \;=\; E\left(\frac{zY}{1-z+zY}\right), \qquad \text{where } 0 \le Y \sim \boldsymbol{Q}$$

$$=\; \int_0^1 \psi_y(z)\,\boldsymbol{Q} + \int_1^{\infty}\{1-\psi_{\frac{1}{y}}(1-z)\}\,\boldsymbol{Q}(dy),$$

where $\psi_y(z) = Ez^{N_y}$, $N_y \sim Geometric(y)$, $0 \le y \le 1$ with mass function $P(N_y = k) = y(1-y)^{k-1}$, $k \ge 1$. Thus

$$(2.3) \qquad \phi(z) \;\equiv\; \int_0^1 Ez^{N_y}\boldsymbol{Q}(dy) + \int_1^{\infty}\{1-\; E(1-z)^{N_{\frac{1}{y}}}\}\boldsymbol{Q}(dy)$$

for $0 < z \le 1$. Hence, for the EMP Shock Model, the generating function of the shock resistance probabilities $\bar{P}_k = P(J > k)$, must be

$$M(z) \equiv \sum_{k=0}^{\infty}\bar{P}_k z^k \;=\; \frac{1-\phi(z)}{1-z}, \qquad 0 < z < 1$$

$$(2.4) \qquad\qquad =\; \int_0^1 \frac{1-Ez^{N_y}}{1-z}\boldsymbol{Q}(dy) + \int_1^{\infty}E(1-z)^{N_y^*}\boldsymbol{Q}(dy),$$

where $N_y^* \stackrel{d}{=} N_{\frac{1}{y}}-1$, $y > 1$.

The integrand in the first term of (2.4) is easily seen to be $\{1-z(1-y)\}^{-1} = \sum_{k=0}^{\infty}(1-y)^k z^k$ for $0 < y \le 1$ and $0 < z < 1$. To evaluate the second term in (2.4), use Lemma 2.2, to get

$$q_n = P(N_y^* = n) = P(N_{\widetilde{\Pi}} = n+1) = \widetilde{\Pi}(1-\widetilde{\Pi})^n.$$

where $\widetilde{\Pi} = y^{-1} < 1$ and $n = 0,1,2,\dots$. The coefficients $c_k$, defined in Lemma 2.2,

in the expansion of $E(1-z)^{N_y^*}$, are

$$
\begin{aligned}
c_k &= (-1)^k \widetilde{\Pi} \sum_{j=0}^{\infty} \binom{k+j}{k}(1-\widetilde{\Pi})^{k+j} \\
&= (-1)^k \frac{\widetilde{\Pi}(1-\widetilde{\Pi})^k}{\widetilde{\Pi}^{k+1}} \left[\sum_{j=0}^{\infty} \binom{k+j}{k} \widetilde{\Pi}^k (1-\widetilde{\Pi})^j \widetilde{\Pi}\right] \\
&= (-1)^k \frac{\widetilde{\Pi}(1-\widetilde{\Pi})^k}{\widetilde{\Pi}^{k+1}} \sum_{j=0}^{\infty} \mathrm{P}((k+1)\text{st success in}(k+j+1)\text{st trial})
\end{aligned}
$$

$$
(2.5) \qquad = (-1)^k \left(\frac{1-\widetilde{\Pi}}{\widetilde{\Pi}}\right)^k = (1-y)^k,
$$

since $\widetilde{\Pi} = y^{-1} < 1$. Note, the right hand side of (2.5) is positive for $k$ even and negative for $k$ odd. Thus (2.4) can be rewritten as,

$$
(2.6) \qquad \sum_{k=0}^{\infty} \bar{P}_k z^k = \int_0^1 \left\{\sum_{k=0}^{\infty}(1-y)^k z^k\right\} \boldsymbol{Q}(dy)
$$
$$
+ \int_1^{\infty} \left\{\sum_{k=0}^{\infty}(-1)^k(y-1)^k z^k\right\} \boldsymbol{Q}(dy).
$$

In the first term of the right hand side of (2.6), the series is absolutely convergent. Hence the integral and the summation can be interchanged. But the series in the second term converges iff $y < 2$. This implies, if $\boldsymbol{Q}[2,\infty) > 0$ then the right hand side of (2.6) diverges. Thus, $\boldsymbol{Q}[2,\infty) = 0$. $\qquad \square$

Theorem 2.4 allows us to have a representation of the shock resistance probabilities $\bar{P}_k$ via the mixing distribution $\boldsymbol{Q}$. From (2.6), note that, for all $z \in (0,1)$,

$$
(2.7) \qquad \sum_{k=0}^{\infty} \bar{P}_k z^k = \sum_{k=0}^{\infty} \left\{\int_0^2 a^k(y)\boldsymbol{Q}(dy)\right\} z^k,
$$

where

$$
(2.8) \qquad a(y) = 1-y = \begin{cases} (1-y) > 0, & \text{for } 0 \le y < 1, \\ -(y-1) < 0, & \text{for } 1 < y < 2. \end{cases}
$$

Equating coefficients of $z^k$ on both sides of (2.7), for $k = 0, 1, 2, \ldots$ we have

$$
(2.9) \qquad \bar{P}_k = \int_0^2 a^k(y)\boldsymbol{Q}(dy), \qquad k = 0, 1, 2, \ldots.
$$

Theorem 2.4 implies that for the function $\phi(z)$ defined by (1.1) to be a p.g.f., we must have $\boldsymbol{Q}[2,\infty)=0$. Is this the sharpest possible result? Or is there a sharper necessary condition? We do not know. A counterexample at the end of this section exhibits a mixing distribution $\boldsymbol{Q}$ with support strictly larger than the unit interval, but the right end point of support also falls strictly short of 2. The counterexample does establish however that, as expected, the full converse to Theorem 2.4 is false.

For what conditions on the mixing distribution $\boldsymbol{Q}$, is $\bar{P}_k$ a tail probability of a discrete non-negative random variable (not necessarily discrete SDFR)? In other words, what conditions on $\boldsymbol{Q}$ would ensure that $\bar{P}_k$ is non-negative and monotonically decreasing? Note, we already have $\bar{P}_k \to 0$ directly from (2.9). Interestingly, a weak converse to Theorem 2.4 does hold under a condition which turns out to be directly relevant to our original question posed in Section 1.

**Theorem 2.5.** *Suppose that the shock model distribution $S$ in (2.1) is SDFR. Then the corresponding number of shocks to failure has the discrete SDFR property if and only if the support of the mixing distribution $\boldsymbol{Q}$ in the representation of its p.g.f. in (2.2) is contained in $(0, 1]$.*

*Proof.* Given that the shock model probability $S$ in (2.1) is nondiscrete SDFR, we have the representation in (2.9) of the shock resistance probabilities, which can be expressed as

$$(2.10) \qquad \bar{P}_k = E\{(1-Y)^k \, 1_{\{0 \leq Y < 1\}} + (-1)^k (Y-1)^k \, 1_{\{1 \leq Y < 2\}}\} = EV^k,$$

where $V = 1 - Y$ and $Y \sim \boldsymbol{Q}$ with support in $[0, 2)$. Thus $-1 < V < 1$, w.p. 1, with d.f. $\boldsymbol{F}$, given by

$$\boldsymbol{F}(v) = 1 - \boldsymbol{Q}(1 - v - 0), \quad \text{for } -1 < v < 1.$$

Accordingly, $\bar{P}_k$ is a moment sequence

$$(2.11) \qquad \bar{P}_k = \int_{\{|v|<1\}} v^k \boldsymbol{F}(dv), \quad k = 0, 1, 2, \dots,$$

of a distribution $\boldsymbol{F}$ with support in $(-1, 1)$. On the other hand, by the classic result of Hausdorff (Feller [4]), $\bar{P}_k$ is discrete SDFR (or, a completely monotone sequence satisfying $\bar{P}_0 = 1$) $\iff \bar{P}_k$ is the moment sequence of a distribution with support in $[0, 1]$. Since d.f.s with bounded support are uniquely determined by their moment sequence, this implies $V$ must have support in $[0, 1]$. Hence,

$$0 = P(-1 < V < 0) = P(1 < Y < 2).$$

In conjunction with Theorem 2.4, this implies that $\boldsymbol{Q}$ must be supported by the unit interval.

Conversely, suppose $\boldsymbol{Q}(0, 1] = 1$. Then the representation (2.9) of the shock resistance probabilities reduces to

$$\bar{P}_k = \int_0^1 (1-y)^k \, \boldsymbol{Q}(dy), \quad k = 0, 1, 2, \dots.$$

Thus $\bar{P}_k$ constitutes the moment sequence of a unique probability measure in the unit interval. Again, in virtue of Hausdoff's theorem, as remarked in the proof of Theorem 2.2, this is equivalent to $\bar{P}_k$ being discrete SDFR. □

Contrary to crude intuition, the mixing distribution $\boldsymbol{Q}$ can have positive mass in $(1, 2)$ and the corresponding $\phi(z)$ can still be a p.g.f. as Theorem 2.4 suggests and the following counterexample demonstrates. In view of Theorem 2.5 however, it can no longer have the discrete SDFR property.

**Counterexample.** We construct a mixing distribution $\boldsymbol{Q}$ on the half line such that (2.9) yields a tail probability sequence. Note that in view of (2.10), it is enough to find a distribution $\boldsymbol{F}$ of a random variable $V$ with support in $(-1, 1)$ such that its moment sequence $EV^k$ is monotonically decreasing in $k = 0, 1, 2, \ldots$, which requires

$$EV^{2n} \geq EV^{2n+1} \geq EV^{2n+2}, \quad n = 0, 1, 2, \ldots, \text{ where } -1 < V < 1 \text{ w.p. } 1.$$

Only the second inequality needs to be verified for a counterexample, since the first one is free. Setting $U = -V$, we see that $EV^{2n} \geq EV^{2n+1}$ iff

$$\begin{aligned}
&E\{1_{\{0 < U < 1\}} U^{2n} + 1_{\{0 \leq V < 1\}} V^{2n}\} \\
&\geq E\{-1_{\{0 < U < 1\}} U^{2n+1} + 1_{\{0 \leq V < 1\}} V^{2n+1}\} \\
&\Leftrightarrow E\{1_{\{0 < U < 1\}} U^{2n}(1 + U)\} \\
(2.12) \qquad &\geq -E\{1_{\{0 \leq V < 1\}} V^{2n}(1 - V)\},
\end{aligned}$$

which clearly holds, since the left hand side is positive, while the right hand side is non-positive. Now set $V$ to have a density given by

$$(2.13) \qquad f(v) = \frac{\beta}{\alpha} 1_{\{-\alpha \leq v < 0\}} + (1 - \beta) 1_{\{0 \leq v < 1\}},$$

with discontinuities at $v = -\alpha$ and $0$. The corresponding d.f. with parameters $\alpha > 0, \beta > 0$, is

$$\boldsymbol{F}(v) = \begin{cases}
0, & \text{if } v < -\alpha, \\
\beta(1 + v\alpha^{-1}), & \text{if } -\alpha \leq v < 0, \\
\beta + (1 - \beta)v, & \text{if } 0 \leq v < 1, \\
1, & \text{if } v \geq 1.
\end{cases}$$

We need to verify that there exist a choice of $(\alpha, \beta) \in (0, 1)^2$ such that $EV^{2n+1} \geq EV^{2n+2}$ for every $n \geq 0$. Computations analogous to those leading to (2.12) shows that this is equivalent to demanding

$$\int_0^1 (v^{2n+1} + v^{2n+2}) f(-v) dv \leq \int_0^1 (v^{2n+1} - v^{2n+2}) f(v) dv.$$

For the density in (2.13), the left hand side equals

$$\frac{\beta}{\alpha} \int_0^\alpha (v^{2n+1} + v^{2n+2}) \, dv = \beta \alpha^{2n+1} \frac{2n(1 + \alpha) + (3 + 2\alpha)}{(2n + 2)(2n + 3)},$$

while the right hand side equals $(1 - \beta)\{(2n + 2)(2n + 3)\}^{-1}$. Thus we need to show that it is possible to have,

$$(2.14) \qquad \beta \alpha^{2n+1} \{2n(1 + \alpha) + (3 + 2\alpha)\} \leq (1 - \beta), \quad n = 0, 1, 2, \ldots.$$

If we choose $(\alpha, \beta)$ such that

$$(2.15) \qquad \alpha\beta \leq \frac{(1 - \beta)}{(3 + 2\alpha)},$$

then (2.14) would be assured under (2.15), provided

$$(2.16) \qquad 2n \, \alpha^{2n}(1 + \alpha) \leq \frac{(1 - \beta)}{\alpha\beta}(1 - \alpha^{2n}), \quad n \geq 1.$$

This is indeed feasible for a wide range of $(\alpha, \beta) \in (0,1)^2$. For example take $1/3 \leq \beta < 1$, and choose $\alpha > 0$ such that $\alpha(3 + 2\alpha) < 2$, thus satisfying (2.15). In particular, choosing $0 < \alpha < 2/7$ is enough, since this would imply

$$0 < \alpha < 2/7 < \frac{2}{3 + 2\alpha}.$$

For such a choice of $\alpha$, the right hand side of (2.16) is bounded below by $\{7(1 - (2/7)^{2n}\}$, whereas the left hand side is bounded above by

$$2n(2/7^{2n})(9/7) < 3n(2/7)^{2n} < 3n(1/3)^{2n} = n3^{-(2n-1)} < 1 < \{7(1 - (2/7)^{2n})\},$$

the lower bound to the right hand side of (2.16). This establishes our claim that the class of d.f.s with support in $(-1, 1)$ whose moment sequence defines a discrete tail probability is a nonempty class.

Finally, to get a counterexample of a mixing distribution $\boldsymbol{Q}$ with support in $[0, 2)$ such that $\boldsymbol{Q}(1, 2) > 0$, take $\boldsymbol{Q}$ to be the d.f. of $Y = 1 - V$ where $V$ has density as in (2.13). The corresponding continuous d.f. is $\boldsymbol{Q}(y) = P(V \geq 1 - y) = 1 - \boldsymbol{F}(1 - y)$, for $0 < y < 1 + \alpha$. Thus;

$$\boldsymbol{Q}(y) = \begin{cases} 0, & \text{if } y < 0, \\ (1 - \beta)y, & \text{if } 0 \leq y < 1, \\ (1 - \beta) + \alpha^{-1}\beta(y - 1), & \text{if } 1 \leq y < (1 + \alpha), \\ 1, & \text{if } y \geq 1 + \alpha. \end{cases}$$

Using this in (2.9) leads to,

$$\bar{P}_k = \frac{1 - \beta}{k + 1} + (-1)^k \frac{\beta\alpha^k}{k + 1}, \quad k = 0, 1, 2, \dots$$

and defines a tail probability for $0 < \alpha < 2/7$ and $\beta > 1/3$, as we have argued. Note, the support of $\boldsymbol{Q}$ exceeds the unit interval. In particular, choosing $\beta = 2/3$ and $\alpha = 1/7$, we have $\Delta^2\bar{P}_1 = -0.00523$. Thus $\bar{P}_k$ is not a completely monotone sequence, and hence is *not* SDFR.

In closing, we note that when the integral in (2.2) is indeed a p.g.f., it can be bounded below (and above, when $EY \leq 1$) by simple geometric p.g.f.s. The Laplace-Stieltjes transform of the EMP shock model distribution $S$ in (2.1) can be correspondingly bounded by Laplace transforms of suitable exponential distributions. In stochastic comparison parlance, these findings can be stated in terms of stochastic domination via the "*generating function order*" $<_g$ and "*Laplace transform order*" $<_L$ (for definitions, see Müller and Stoyan [5]).

Let $N_\alpha$ be a geometric r.v. on the positive integers, with mass function $P(N_\alpha = k) = \alpha(1 - \alpha)^{k-1}$, $k \geq 1$, as defined in proof of Theorem 2.4, and let $(exp)_\mu$ denote an exponential distribution with mean $\mu$.

**Theorem 2.6.** *For the EMP shock model distribution $S$ in (2.1) and the corresponding number of shocks $(J)$ to failure, we have*

$$J <_g N_{(EJ)^{-1}} \quad and \quad S <_L (exp)_{\mu = \lambda^{-1}EJ}$$
$$If \quad EY \leq 1, \text{ then, } N_{(EY)} <_g J \quad and \quad (exp)_{\mu = (\lambda EY)^{-1}} <_L S.$$

*The bounds are sharp, with equality attained when $J = 1$ w.p. 1.*

*Proof.* Note that the expected number of shocks to failure, satisfies

$$
\begin{aligned}
1 \leq EJ & = \lim_{z \to 1-} \frac{1 - \phi(z)}{1 - z} \\
(2.17) \qquad & = \lim_{z \to 1-} \int_0^\infty \frac{\boldsymbol{Q}(dy)}{1 - z + zy} = \int_0^\infty \frac{\boldsymbol{Q}(dy)}{y},
\end{aligned}
$$

by the dominated convergence theorem, since the integrand is $\downarrow$ in $y$, for each $z \in (0,1)$ and is thus bounded above by the constant $(1-z)^{-1}$. We note in passing that, using the necessary condition in Theorem 2.4, we could argue : $EJ = \int_0^2 y^{-1} \boldsymbol{Q}(dy) \geq \frac{1}{2}$, which however is superseded by the sharper bound $EJ \geq 1$, since $J \geq 1$, w.p. 1. In fact (2.17) together with Theorem 2.4 now implies $\int_0^2 y^{-1} \boldsymbol{Q}(dy) \geq 1$, which implicitly restricts the amount of mass $\boldsymbol{Q}$ can allocate to the interval $(1,2)$.

Finally, appealing again to the concavity in $y$, and hence convexity in $x \equiv y^{-1}$, of the integrand in (2.2), for each $z \in (0,1)$; standard Jensen's inequalities give,

$$
(2.18) \qquad \frac{z(EY)}{1 - z + z(EY)} \geq \phi(z) = E\left(\frac{z}{z + (1-z)Y^{-1}}\right) \geq \frac{z}{z + (1-z)(EJ)},
$$

using (2.17). We note that while the upper bound is always valid, it is a geometric p.g.f. only when $EY \leq 1$. The lower bound is of course the p.g.f. of a geometric distribution with parameter $(EJ)^{-1}$.

The Laplace ordering claims now follow from (2.18), by using Lemma 2.1 and the substitution $z = \lambda/(\lambda + s)$. All bounds are easily seen to be sharp, being attained when the first shock causes failure (when $J = 1$, w.p. 1). $\qquad \square$

**Acknowledgment.** It is a pleasure to dedicate this paper in honor of Professor P. K. Sen; whose constructive comments, while the work was in progress, are gratefully acknowledged.

## References

[1] Bhattacharjee, M. C. (2005). Strong versions of the DFR property. CAMS Research Report 0405-18. Available at http://m.njit.edu/CAMS/Technical_Reports/CAMS04_05/report18.pdf.

[2] Esary, J. D. and Marshall, A. W. (1973). Shock models and wear processes. *Ann. Probab.* **1** 624–649. MR0350893

[3] Feller, W. (1939). Completely montone functions and sequences. *Duke J. Math.* **5** 661–674. MR0000315

[4] Feller, W. (1960). *An Introduction to Probability Theory and Its Application.* II. Wiley, New York.

[5] Müller, A. and Stoyan, D. (2002). *Comparison Methods for Stochastic Models and Risk.* Wiley, New York. MR1889865

[6] Pakes, A. G. (1997). On the recognition and structure of probabilitygenerating functions. In *Classical and Modern Branching Processes* (K. B. Athreya and P. Jagers, eds.) 263–284. Springer, New York. MR1601686